\documentclass{article}
\usepackage{amsfonts}
\usepackage{amsmath}
\usepackage{color}

\newtheorem{theorem}{Theorem}

\newtheorem{example}[theorem]{Example}

\newtheorem{lemma}[theorem]{Lemma}

\begin{document}

\title{Branching processes in random environment with immigration stopped at
zero\thanks{%
\noindent Doudou Li and Mei Zhang were supported by the Natural Science
Foundation of China under the grant 11871103, V.Vatutin was supported by the
High-End Foreign Experts Recruitment Program (No. GDW20171100029).}}
\author{Elena Dyakonova\thanks{%
Steklov Mathematical Institute, 8 Gubkin St., Moscow, 119991, Russia Email:
elena@mi-ras.ru}, Doudou Li\thanks{%
School of Mathematical Sciences \& Laboratory of Mathematics and Complex
Systems, Beijing Normal University, Beijing 100875, P.R. China. Email:
lidoudou@mail.bnu.edu.cn},~Vladimir Vatutin\thanks{%
Steklov Mathematical Institute, 8 Gubkin St., Moscow, 119991, Russia and
Beijing Normal University, Beijing 100875, P.R. China Email:
vatutin@mi-ras.ru} \ and Mei Zhang\thanks{%
School of Mathematical Sciences \& Laboratory of Mathematics and Complex
Systems, Beijing Normal University, Beijing 100875, P.R. China. Email:
meizhang@bnu.edu.cn} }
\date{}
\maketitle

\begin{abstract}
A critical branching process with immigration which evolve in a random
environment is considered. Assuming that immigration is not allowed when
there are no individuals in the aboriginal population we investigate the
tail distribution of the so-called life period of the process, i.e., the
length of the time interval between the moment when the process is initiated
by a positive number of particles and the moment when there are no
individuals in the population for the first time.
\end{abstract}


\section{~ Introduction and statement of main results}

\setcounter{theorem}{0}

We consider branching processes allowing immigration and evolving in a
random environment. In such a process individuals reproduce independently of
each other according to random offspring distributions which vary from one
generation to the other. In addition, immigrants arrive to each generation
independently on the development of the population and according to the laws
varying at random from generation to generation. To give a formal definition
let $\Delta =\left( \Delta _{1},\Delta _{2}\right) $ be the space of all
pairs of probability measures on $\mathbb{N}_{0}=\{0,1,2,\ldots \}.$
Equipped with the component-wise metric of total variation $\Delta $ becomes
a Polish space. Let $\mathbf{Q}=\{F,G\}$ be a random vector with
independent components taking values in $\Delta $, and let $\mathbf{Q}%
_{n}=\{F_{n},G_{n}\},n=1,2,\ldots ,$ be a sequence of independent copies of $%
\mathbf{Q}$. The infinite sequence $\mathcal{E}=\left\{ \mathbf{Q}_{1},%
\mathbf{Q}_{2},...\right\} $ is called a random environment.

A sequence of $\mathbb{N}_{0}$-valued random variables $\mathbf{Y}=\left\{
Y_{n},\ n\in \mathbb{N}_{0}\right\} $ specified on\ the respective
probability space $(\Omega ,\mathcal{F},\mathbf{P})$ is called a branching
process with immigration in the random environment (BPIRE), if $Y_{0}$ is
independent of $\mathcal{E}$ and, given $\mathcal{E}$ the process $\mathbf{Y}
$ is a Markov chain with
\begin{equation}
\mathcal{L}\left( Y_{n}|Y_{n-1}=y_{n-1},\mathcal{E}=(\mathbf{q}_{1},\mathbf{q%
}_{2},...)\right) =\mathcal{L}(\xi _{n1}+\ldots +\xi _{ny_{n-1}}+\eta _{n})
\label{BasicDefBPimmigr}
\end{equation}%
for every $n\in \mathbb{N}:=\mathbb{N}_{0}\backslash \left\{ 0\right\} $, $%
y_{n-1}\in \mathbb{N}_{0}$ and $\mathbf{q}_{1}=\left( f_{1},g_{1}\right) ,%
\mathbf{q}_{2}=\left( f_{2},g_{2}\right) ,...\in \mathbf{Q}$, where $\xi
_{n1},\xi _{n2},\ldots $ are i.i.d. random variables with distribution $%
f_{n} $ and independent of the random variable $\eta _{n}$ with distribution
$g_{n} $. In the language of branching processes $Y_{n-1}$ is the $(n-1)$th
generation size of the population, $f_{n}$ is the distribution of the number
of children of an individual at generation $n-1$ and $g_{n}$ is the
reproduction law of immigrants at generation $n$.

Along with the process $\mathbf{Y}$ we consider a branching process $\mathbf{%
Z}=\left\{ Z_{n},\ n\in \mathbb{N}_{0}\right\} $ in the random environment $%
\mathcal{E}_{1}=\left\{ F_{1},F_{2},...\right\} $ which, given $\mathcal{E}%
_{1}$ is a Markov chain with $Z_{0}=1$ and, for $n\in \mathbb{N}$
\begin{equation}
\mathcal{L}\left( Z_{n}|Z_{n-1}=z_{n-1},\mathcal{E}_{1}=(f_{1},f_{2},...)%
\right) =\mathcal{L}(\xi _{n1}+\ldots +\xi _{nz_{n-1}}).  \label{BPordinary}
\end{equation}

It will be convenient to assume that if $Y_{n-1}=y_{n-1}>0$ is the
population size of the ($n-1)$th generation of $\mathbf{Y}$ then first $\xi
_{n1}+\ldots +\xi _{ny_{n-1}}$ individuals of the $n$th generation are born
and than $\eta _{n}$ immigrants enter the population.

This agreement allows us to consider a modified version $\mathbf{W}=\left\{
W_{n},\ n\in \mathbb{N}_{0}\right\} $ of the process $\mathbf{Y}$ specified
as follows. Assume, without loss of generality that $Y_{0}>0.$ Let $%
W_{0}=Y_{0}$ and for $n\geq 1$,
\begin{equation}
W_{n}:=\left\{
\begin{array}{cc}
0, & \text{ if }T_{n}:=\xi _{n1}+\ldots +\xi _{nW_{n-1}}=0, \\[0.1in]
T_{n}+\eta _{n}, & \text{if }T_{n}>0.%
\end{array}%
\right. .  \label{BPstopped}
\end{equation}%
We call $\mathbf{W}$ as a branching process with immigration stopped at zero
and evolving in the random environment.

The aim of the present paper is to study the tail distribution of the random
variable
\begin{equation*}
\zeta :=\min \left\{ n\geq 1:W_{n}=0\right\}
\end{equation*}%
under the annealed approach. To formulate our main result we consider the
so-called associated random walk $\mathbf{S}=\left( S_{0},S_{1},...\right) $%
. This random walk has initial state $S_{0}$ and increments $%
X_{n}=S_{n}-S_{n-1}$, $n\geq 1$, defined as
\begin{equation*}
X_{n}:=\log \mathfrak{m}\left( F_{n}\right)
\end{equation*}%
which are i.i.d. copies of the logarithmic mean offspring number $X:=\log $ $%
\mathfrak{m}(F)$ with%
\begin{equation*}
\mathfrak{m}(F):=\sum_{j=0}^{\infty }jF\left( \left\{ j\right\} \right) .
\end{equation*}%
We suppose that $X$ is a.s. finite.

With each pair of measures $(F,G)$ we associate the respective probability
generation functions
\begin{equation*}
F(s):=\sum_{j=0}^{\infty }F\left( \left\{ j\right\} \right) s^{j},\qquad
G(s):=\sum_{j=0}^{\infty }G\left( \left\{ j\right\} \right) s^{j}.
\end{equation*}%
We impose the following restrictions on the distributions of $F$ and $G$.

\textbf{Hypothesis A1}. The probability generating function $F(s)$ is
geometric with probability 1, that is%
\begin{equation}
F(s)=\frac{q}{1-ps}=\frac{1}{1+\mathfrak{m}(F)(1-s)}  \label{Frac_generating}
\end{equation}%
with random $p,q\in (0,1)$ satisfying $p+q=1$ and
\begin{equation*}
\mathfrak{m}(F)=\frac{p}{q}=e^{\log (p/q)}=e^{X}.
\end{equation*}

\textbf{Hypothesis A2}. There exist real numbers $\kappa \in \lbrack 0,1)$
and $\gamma ,\sigma \in (0,1]$ such that, with probability 1

1) the inequality $F(0)\geq \kappa $ is valid;

2) the estimate%
\begin{equation}
G(s)\leq s^{\gamma }  \label{Ineq1}
\end{equation}%
holds for all $s\in \lbrack \kappa ^{\sigma },1].$

To formulate one more assumption we set
\begin{equation*}
M_{n}:=\max \left( S_{1},...,S_{n}\right) ,\quad L_{n}:=\min \left(
S_{0},S_{1},...,S_{n}\right) ,
\end{equation*}%
and, given $S_{0}=0,$ introduce the right-continuous function $U:\mathbb{R}$
$\rightarrow [0,\infty )$ specified by the relation%
\begin{equation}
U(x):=I\left\{ x\geq 0\right\} +\sum_{n=1}^{\infty }\mathbf{P}\left(
S_{n}\geq -x,M_{n}<0\right) ,  \label{DefU}
\end{equation}%
where $I(A)$ is the indicator of the event $A.$

One may check (see, for instance, \cite{4h}\ and \cite{ABKV}) that for any
oscillating random walk
\begin{equation}
\mathbf{E}\left[ U(x+X);X+x\geq 0\right] =U(x),\quad x\geq 0.  \label{Mes1}
\end{equation}

\textbf{Hypothesis} $\mathbf{A3}$. The distribution of $X$ is nonlattice,
the sequence $\left\{ S_{n},n\geq 0\right\} $ satisfies the Doney-Spitzer
condition
\begin{equation}
\lim_{n\rightarrow \infty }\mathbf{P}\left( S_{n}>0\right) =:\rho \in (0,1),
\label{Don1}
\end{equation}%
and there exists $\varepsilon >0$ such that
\begin{equation}
\mathbf{E}\left( \log ^{+}G^{\prime }(1)\right) ^{\rho ^{-1}+\varepsilon
}<\infty \quad \text{ and \ \ }\mathbf{E}\left( U(X)\log ^{+}G^{\prime
}(1)\right) ^{1+\varepsilon }<\infty ,  \label{G_restriction}
\end{equation}%
where $\log ^{+}x=\max \left( 0,\log x\right) $.

We now formulate our main result.

\begin{theorem}
\label{T_infvar0}Let Hypotheses A1 - A3 be satisfied. Then there exists a
function $l(n)$ slowly varying at infinity such that
\begin{equation*}
\mathbf{P}\left( \zeta >n\right) \sim \frac{l(n)}{n^{1-\rho }}
\end{equation*}%
as $n\rightarrow \infty .$
\end{theorem}

It is convenient to describe the range of possible values of the parameter $%
\kappa $ by examples.

Let
\begin{equation*}
\mathcal{A}:=\{0<\alpha <1;\,|\beta |<1\}\cup \{1<\alpha <2;|\beta |\leq
1\}\cup \{\alpha =1,\beta =0\}\mathbf{\cup }\{\alpha =2,\beta =0\}
\end{equation*}%
be a subset in $\mathbb{R}^{2}.$ For $(\alpha ,\beta )\in \mathcal{A}$ and a
random variable $X$ we write $X\in \mathcal{D}\left( \alpha ,\beta \right) $
if the distribution of $X$ belongs to the domain of attraction of a stable
law with characteristic function%
\begin{equation}
\mathcal{G}_{\alpha ,\beta }\mathbb{(}t\mathbb{)}:=\exp \left\{
-c|t|^{\,\alpha }\left( 1-i\beta \frac{t}{|t|}\tan \frac{\pi \alpha }{2}%
\right) \right\} ,\ c>0,  \label{std}
\end{equation}%
and, in addition, $\mathbf{E}\left[ X\right] =0$ if this moment exists. If $%
X_{n}\overset{d}{=}X\in \mathcal{D}\left( \alpha ,\beta \right) $ then the
parameter $\rho $ in (\ref{Don1}) is given by the formula (see, for
instance, \cite{Zol57})
\begin{equation}
\displaystyle\rho =\left\{
\begin{array}{ll}
\frac{1}{2},\ \text{if \ }\alpha =1, &  \\
\frac{1}{2}+\frac{1}{\pi \alpha }\arctan \left( \beta \tan \frac{\pi \alpha
}{2}\right) ,\text{ otherwise}. &
\end{array}%
\right.  \label{ro}
\end{equation}

Note that if $\mathbf{E}\left[ X\right] =0$ and $\mathbf{Var}X\in (0,\infty
) $ then the central limit theorem implies $\rho =1/2$.

\setcounter{theorem}{0}

\begin{example}
If Hypothesis A1 is valid and
\begin{equation*}
X=\log \mathfrak{m}(F)=\log (p/q)\in \mathcal{D}\left( \alpha ,\beta \right)
\end{equation*}%
with $\alpha \in (0,2)$ then
\begin{equation}
\mathbf{P}\left( \log (p/q)>x\right) \sim \frac{1}{x^{\alpha }l_{1}(x)}\quad
\text{as }x\rightarrow \infty ,  \label{Tailtwo}
\end{equation}%
where $l_{1}(x)$ is a function slowly varying at infinity. Therefore,
\begin{equation*}
\mathbf{P}\left( \log \frac{q}{1-q}<-x\right) \sim \frac{1}{x^{\alpha
}l_{1}(x)}
\end{equation*}%
as $\ x\rightarrow \infty $ implying%
\begin{equation*}
\mathbf{P}\left( F(0)=q<\frac{e^{-x}}{1+e^{-x}}\right) \sim \frac{1}{%
x^{\alpha }l_{1}(x)}.
\end{equation*}%
As a result, $\mathbf{P}\left( F(0)<y\right) >0$ for any $y>0.$

Thus, if $\alpha \in (0,2)$ then point 1) of Hypothesis A2 reduces to the
trivial inequality $F(0)\geq \kappa =0$. Moreover, given $\kappa =0$ point
2) of Hypothesis A2 implies $G\left( 0\right) =0$ which, in turn, leads to
the inequality
\begin{equation*}
G(s)=\sum_{j=1}^{\infty }G\left( \left\{ j\right\} \right) s^{j}\leq s
\end{equation*}%
for all $s\in \lbrack 0,1]$. The last means that at least one immigrant
enters $\mathbf{W}$ each time when it is allowed by (\ref{BPstopped}).
\end{example}

The case $\mathbf{E}\left[ X^{2}\right] <\infty $ is less restrictive and
allows for \ $\kappa >0$, i.e., for the absence of immigrants in some
generations of $\mathbf{W}$ (even they are allowed).

\begin{example}
Let
\begin{equation*}
F(s)=\left\{
\begin{array}{ccc}
\frac{1}{1+63\left( 1-s\right) } & \text{with probability} & \frac{1}{2}, \\
&  &  \\
\frac{63}{64-s} & \text{with probability} & \frac{1}{2}%
\end{array}%
\right.
\end{equation*}%
and the probability generating function of immigrants be deterministic:
\begin{equation*}
G(s)=\frac{2}{3}s^{2}+\frac{1}{3}\text{ with probability 1.}
\end{equation*}%
Clearly, $\mathbf{E}\left[ \log \mathfrak{m}(F)\right] =0$, $\mathbf{Var}%
\left[ \log \mathfrak{m}(F)\right] \in \left( 0,\infty \right) $. It is not
difficult to see that
\begin{equation*}
F(0)\geq 1/64\text{ and \ }G(s)\leq s^{1/3}\text{ for all }s\in \left[
8^{-1},1\right] =\left[ 64^{-1/2},1\right] .
\end{equation*}%
Thus, the conditions of Theorem \ref{T_infvar0} fulfill with $\kappa =1/64$,
$\gamma =1/3$ and $\sigma =1/2.$
\end{example}

We note that Zubkov \cite{Zub72} considered a problem similar to ours for a
branching process with immigration $\left\{ Y_{c}(n),n\geq 0\right\} $
evolving in a constant environment. He assumed that $G\left( 0\right) >0$
and investigated the distribution of the so-called life period $\zeta _{c}$
of such a process initiated at time $N$ and defined as
\begin{equation*}
Y_{c}(N-1)=0,\min_{N\leq k<N+\zeta _{c}}Y_{c}(k)>0,Y_{c}(N+\zeta _{c})=0.
\end{equation*}%
The same problem for other models of branching processes with immigration
evolving in a constant environment was analysed, for instance, in \cite%
{BM1983}, \cite{Mit82}, \cite{ST83} and \cite{Vat77}.

Various properties of BPIRE were investigated by several authors (see, for
instance, \cite{AFa2013}, \cite{KP74}, \cite{KKS75},\cite{Key87},\cite%
{Rot2007}, \cite{Tan81}). However, asymptotic properties of the life periods
of BPIRE were not considered up to now.

\section{\protect\bigskip Auxiliary statements}

\setcounter{theorem}{0}

Given the environment $\mathcal{E}=\left\{ (F_{n},G_{n}),n\in \mathbb{N}%
\right\} $, we construct the i.i.d. sequence of pairs of generating
functions
\begin{equation*}
F_{n}(s):=\sum_{j=0}^{\infty }F_{n}\left( \left\{ j\right\} \right)
s^{j},\qquad G_{n}(s):=\sum_{j=0}^{\infty }G_{n}\left( \left\{ j\right\}
\right) s^{j}\quad s\in \lbrack 0,1],
\end{equation*}%
and use below the convolutions of the generating functions $F_{1},...,F_{n}$
specified for $0\leq i\leq n-1$ by the equalities%
\begin{eqnarray*}
F_{i,n}(s):= &&F_{i+1}(F_{i+2}(\ldots (F_{n}(s))\ldots )),\quad  \\
F_{n,i}(s):= &&F_{n}(F_{n-1}(\ldots (F_{i+1}(s))\ldots ))\ \text{ and }%
F_{n,n}(s):=s.
\end{eqnarray*}

The evolution of the BPIRE defined by (\ref{BPstopped}) may be now described
for $n\geq 1$ by the relation
\begin{eqnarray}
\mathbf{E}[s^{W_{n}}|\mathcal{E},W_{n-1}] &=&(F_{n}(0))^{W_{n-1}}+\left(
(F_{n}(s))^{W_{n-1}}-(F_{n}(0))^{W_{n-1}}\right) G_{n}(s)  \notag \\
&=&(F_{n}(0))^{W_{n-1}}(1-G_{n}(s))+(F_{n}(s))^{W_{n-1}}G_{n}(s)\ .
\label{DefStop1}
\end{eqnarray}

We assume for convenience that $W_{0}=Y_{0}>0$ has the (random) probability
generating function%
\begin{equation*}
N(0;s):=\frac{G_{0}(s)-G_{0}(0)}{1-G_{0}(0)}
\end{equation*}%
where $G_{0}(s)\overset{d}{=}G(s)$. Other classes of the initial
distribution may be considered in a similar way.

Setting
\begin{equation*}
N(n;s):=\mathbf{E}[s^{W_{n}}|\mathcal{E}],\ n\geq 1
\end{equation*}%
we have by (\ref{DefStop1})
\begin{eqnarray}
N(n;s) &=&\mathbf{E}\left[
(F_{n}(0))^{W_{n-1}}(1-G_{n}(s))+(F_{n}(s))^{W_{n-1}}G_{n}(s)|\mathcal{E}%
\right]  \notag \\
&=&N(n-1;F_{n}(0))(1-G_{n}(s))+N(n-1;F_{n}(s))G_{n}(s)  \label{DefExtin} \\
&=&N(n-1;F_{n}(0))(1-G_{n}(s))+N(n-2;F_{n-1}(0))(1-G_{n-1}(F_{n}(s)))G_{n}(s)
\notag \\
&&+N(n-2;F_{n-1}(F_{n}(s)))G_{n-1}(F_{n}(s))G_{n}(s),  \notag
\end{eqnarray}%
where for $n=1$ one should take into account only the first two equalities.
Assuming $\prod_{j=n+1}^{n}=1$ we obtain by induction
\begin{eqnarray*}
N(n;s)
&=&\sum_{k=0}^{n-1}N(n-k-1;F_{n-k}(0))(1-G_{n-k}(F_{n-k,n}(s)))%
\prod_{j=n-k+1}^{n}G_{j}(F_{j,n}(s)) \\
&&+N(0;F_{0,n}(s))\prod_{j=1}^{n}G_{j}(F_{j,n}(s)).
\end{eqnarray*}%
Note that according to (\ref{DefExtin})
\begin{equation*}
N(n;0)=N(n-1;F_{n}(0)),\,n\geq 1.
\end{equation*}%
Besides,%
\begin{equation*}
\mathbf{E}N(n;0)=\mathbf{P}\left( W_{n}=0\right) =\mathbf{P}\left( \zeta
\leq n\right) .
\end{equation*}

Hence, setting $s=F_{n+1}(0)$, taking the expectation with respect to the
environment and using the independency of the elements of the environment we
get%
\begin{eqnarray}
\mathbf{E}\left[ N(n+1;0)\right] &=&\sum_{k=0}^{n-1}\mathbf{E}\left[ N(n-k;0)%
\right] \mathbf{E}\left[ \left( 1-G_{n-k}(F_{n-k,n+1}(0)\right)
)\prod_{j=n-k+1}^{n}G_{j}(F_{j,n+1}(0))\right]  \notag \\
&&+\mathbf{E}\left[ N(0;F_{0,n+1}(0))\prod_{j=1}^{n}G_{j}(F_{j,n+1}(0))%
\right] .  \label{SubRenewal}
\end{eqnarray}%
Denoting \ for $n\geq 0$
\begin{eqnarray*}
R_{n}:= &&1-\mathbf{E}\left[ N(n;0)\right] =\mathbf{E}\left[ 1-N(n;0)\right]
=\mathbf{P}\left( \zeta >n\right) , \\
H_{n}^{\ast }:= &&\mathbf{E}\left[ \frac{1-G_{0}(F_{0,n+1}(0))}{1-G_{0}(0)}%
\prod_{i=1}^{n}G_{i}(F_{i,n+1}(0))\right] , \\
d_{n}:= &&\mathbf{E}\left[ \prod_{i=1}^{n}G_{i}\left( F_{i,n+1}(0)\right) %
\right] =\mathbf{E}\left[ \prod_{i=1}^{n}G_{i}\left( F_{i,0}(0)\right) %
\right] ,
\end{eqnarray*}%
observing that
\begin{eqnarray*}
H_{n} &:&=\mathbf{E}\left[ (1-G_{0}(F_{0,n+1}(0)))\prod_{i=1}^{n}G_{i}\left(
F_{i,n+1}(0)\right) \right] \\
&=&\mathbf{E}\left[ \prod_{i=1}^{n}G_{i}\left( F_{i,n+1}(0)\right) \right] -%
\mathbf{E}\left[ \prod_{i=1}^{n+1}G_{i}\left( F_{i,n+2}(0)\right) \right]
=d_{n}-d_{n+1},\
\end{eqnarray*}%
and using the equality%
\begin{equation*}
\mathbf{E}\left[ \left( 1-G_{n-k}(F_{n-k,n+1}(0)\right)
)\prod_{j=n-k+1}^{n}G_{j}(F_{j,n+1}(0))\right] =\mathbf{E}\left[ \left(
1-G_{0}(F_{0,k+1}(0)\right) )\prod_{j=1}^{k}G_{j}(F_{j,k+1}(0))\right]
\end{equation*}%
we rewrite (\ref{SubRenewal}) as a renewal type equation
\begin{equation}
R_{n+1}=\sum_{k=0}^{n-1}H_{k}R_{n-k}+H_{n}^{\ast },\ n\geq 0.
\label{Renewal0}
\end{equation}

Let
\begin{equation*}
\mathcal{R}(s):=\sum_{n=1}^{\infty }R_{n}s^{n}.
\end{equation*}

\begin{lemma}
\label{L_generfunct}%
\begin{equation}
\mathcal{R}(s)=\frac{s\mathcal{H}^{\ast }(s)+sR_{1}}{\left( 1-s\right)
D\left( s\right) }  \label{ExplisitR}
\end{equation}%
where%
\begin{equation*}
D\left( s\right) :=\sum_{n=0}^{\infty }d_{n}s^{n}\text{ and }\mathcal{H}%
^{\ast }(s):=\sum_{n=1}^{\infty }H_{n}^{\ast }s^{n}.
\end{equation*}
\end{lemma}

\textbf{Proof}. Set
\begin{equation*}
\mathcal{H}(s):=\sum_{n=0}^{\infty }H_{n}s^{n}.
\end{equation*}%
Clearly,%
\begin{equation*}
s\mathcal{H}(s)=\sum_{n=0}^{\infty }(d_{n}-d_{n+1})s^{n+1}=sD\left( s\right)
-D(s)+1.
\end{equation*}%
Multiplying (\ref{Renewal0}) by $s^{n+1}$ and summing over $n$ from $1$ to $%
\infty $ we get%
\begin{equation*}
\mathcal{R}(s)-sR_{1}=s\mathcal{H}(s)\mathcal{R}(s)+s\mathcal{H}^{\ast }(s)
\end{equation*}%
or%
\begin{equation*}
\mathcal{R}(s)=\frac{s\mathcal{H}^{\ast }(s)+sR_{1}}{1-s\mathcal{H}(s)}=%
\frac{s\left( \mathcal{H}^{\ast }(s)+R_{1}\right) }{\left( 1-s\right)
D\left( s\right) }.
\end{equation*}

The lemma is proved.

Denote for $0\leq i\leq n$
\begin{equation*}
A_{n}:=e^{S_{n}},\quad B_{i,n}:=\sum_{k=i}^{n}e^{S_{k}},\quad B_{n}:=B_{0,n},
\end{equation*}%
and introduce the function
\begin{equation*}
C_{n}(s):=\prod_{i=1}^{n}F_{i,0}(s).
\end{equation*}

\begin{lemma}
\label{L_represent}Under Hypothesis $A1$%
\begin{equation*}
C_{n}:=C_{n}(0)=\frac{1}{B_{n}}.
\end{equation*}
\end{lemma}

\textbf{Proof}. Hypothesis A1 implies
\begin{equation}
\ F_{i}(s)=\frac{q_{i}}{1-p_{i}s}=\frac{1}{1+e^{X_{i}}\left( 1-s\right) }
\label{Frac2}
\end{equation}%
for all $i=1,2,\ldots $. Using these equalities it is not difficult to check
by induction that, for $n\geq 1$
\begin{equation*}
F_{n,0}(s)=1-\frac{A_{n}}{\left( 1-s\right) ^{-1}+B_{1,n}}=\frac{\left(
1-s\right) ^{-1}+B_{1,n-1}}{\left( 1-s\right) ^{-1}+B_{1,n}},
\end{equation*}%
where $B_{1,0}=0$ by definition. Therefore,
\begin{equation}
C_{n}(s)=\prod_{i=1}^{n}\frac{\left( 1-s\right) ^{-1}+B_{1,i-1}}{\left(
1-s\right) ^{-1}+B_{1,i}}=\frac{\left( 1-s\right) ^{-1}}{\left( 1-s\right)
^{-1}+B_{1,n}}.  \label{Repr_gnBar}
\end{equation}

Setting $s=0$ in (\ref{Repr_gnBar}) we prove the lemma.

To go further we need more notation. Let $\mathcal{E}=\left\{ \mathbf{Q}_{1},%
\mathbf{Q}_{2},...\right\} $ be a random environment and let $\mathcal{F}%
_{n},n\geq 1,$ be the $\sigma $-field of events generated by the random
pairs $\mathbf{Q}_{1}=\{F_{1},G_{1}\},\mathbf{Q}_{2}=\{F_{2},G_{2}\},...,%
\mathbf{Q}_{n}=\{F_{n},G_{n}\}$ and the sequence $W_{0},W_{1},...,W_{n}$.
These $\sigma $-fields form a filtration $\mathfrak{F}$. Now the increments $%
\left\{ X_{n},n\geq 1\right\} $ of the random walk $S$ are measurable with
respect to the $\sigma $-field $\mathcal{F}_{n}$. Using the martingale
property (\ref{Mes1}) of $U$ we introduce a sequence of probability measures
$\left\{ \mathbf{P}_{(n)}^{+},n\geq 1\right\} $ on the $\sigma $-field $%
\mathcal{F}_{n}$ by means of the density
\begin{equation*}
d\mathbf{P}_{(n)}^{+}:=U(S_{n})I\left\{ L_{n}\geq 0\right\} d\mathbf{P}.
\end{equation*}%
This and Kolmogorov's extension theorem show that, on a suitable probability
space there exists a probability measure $\mathbf{P}^{+}$ on the $\sigma $%
-field $\mathfrak{F}$ such that (see \cite{4h} and \cite{ABKV} for more
detail)
\begin{equation*}
\mathbf{P}^{+}|\mathcal{F}_{n}=\mathbf{P}_{(n)}^{+},\ n\geq 1.
\end{equation*}

We now formulate two known statements dealing with conditioning $\left\{
L_{n}\geq 0\right\} $.

\begin{lemma}
\label{Lappbasic1}(see Lemma 2.5 in \cite{4h} or Lemma 5.2 in \cite{GV2017})
Let the condition (\ref{Don1}) hold and let $\xi _{1},\xi _{2},\ldots $ be a
sequence of uniformly bounded random variables adapted to the filtration ${%
\mathfrak{F}}$ such that the limit
\begin{equation}
\xi _{\infty }:=\lim_{n\rightarrow \infty }\xi _{n}  \label{AS0}
\end{equation}%
exists $\mathbf{P}^{+}$ - a.s. Then
\begin{equation}
\lim_{n\rightarrow \infty }\mathbf{E}[\xi _{n}\,|\,L_{n}\geq 0]=\mathbf{E}%
^{+}\left[ \xi _{\infty }\right] .  \label{ASSS0}
\end{equation}
\end{lemma}

Let%
\begin{equation*}
\tau (n):=\min \left\{ i\geq 0:S_{i}=L_{n}\right\} .
\end{equation*}

\begin{lemma}
\label{Lappbasic}(see Lemma 2.2 in \cite{4h}) Let $u(x),x\geq 0,$ be a
nonnegative, nonincreasing function with $\int_{0}^{\infty }u(x)dx<\infty $.
If the condition (\ref{Don1}) holds then, for every $\varepsilon >0,$ there
exists a positive number $m=m(\varepsilon )$ such that for all $n\geq m$%
\begin{equation*}
\sum_{k=m}^{n}\mathbf{E}\left[ u(-S_{k});\tau (k)=k\right] \mathbf{P}\left(
L_{n-k}\geq 0\right) \leq \varepsilon \mathbf{P}\left( L_{n}\geq 0\right) .
\end{equation*}
\end{lemma}

\section{Proof of the main result}

It is known (see, for instance, \cite{Rog71} or \cite{BGT87}, Theorem
8.9.12) that if Hypothesis A3~is valid then there exists a slowly varying
function $l_{2}(n)$ such that
\begin{equation}
\mathbf{P}\left( L_{n}\geq 0\right) \sim \frac{l_{2}(n)}{n^{1-\rho }},\quad
n\rightarrow \infty .  \label{Ladder1}
\end{equation}

We now prove an important statement describing the asymptotic behavior of $%
d_{n}$ as $n\rightarrow \infty $. To this aim we introduce the reflected
random walk
\begin{equation*}
\tilde{S}_{0}=0,\ \tilde{S}_{k}=\tilde{X}_{1}+...+\tilde{X}_{k},\ k\geq 1,
\end{equation*}%
where \ $\tilde{X}_{k}=-X_{k}$ and supply in the sequel the relevant
variables and measures by the upper symbol $\symbol{126}$ $.$

Note that $\tilde{X}_{k}\in \mathcal{D}\left( \alpha ,-\beta \right) $ and
\begin{equation*}
\lim_{n\rightarrow \infty }\mathbf{P}\left( \tilde{S}_{n}>0\right)
=\lim_{n\rightarrow \infty }\mathbf{P}\left( S_{n}<0\right) =1-\rho .\text{ }
\end{equation*}%
Hence it follows that
\begin{equation}
\mathbf{P}\left( \tilde{L}_{n}\geq 0\right) \sim \frac{l_{3}(n)}{n^{\rho }}%
,\quad n\rightarrow \infty ,  \label{Ladder2}
\end{equation}%
for a slowly varying function $l_{3}(n)$.

\begin{lemma}
\label{L_infvarFrac}If Hypotheses A1-A3 are satisfied then there exists a
constant $\theta >0$ such that
\begin{equation}
d_{n}\sim \theta \mathbf{P}\left( \tilde{L}_{n}\geq 0\right) \sim \theta
\frac{l_{3}(n)}{n^{\rho }},\quad n\rightarrow \infty .  \label{Asym_zero}
\end{equation}
\end{lemma}

\textbf{Proof}. According to Lemma \ref{L_represent}

\begin{equation*}
C_{n}=\frac{1}{B_{n}}=\frac{1}{1+e^{-\tilde{S}_{1}}+...+e^{-\tilde{S}_{n}}}=:%
\frac{1}{\tilde{B}_{n}}.
\end{equation*}%
We set%
\begin{equation*}
\tilde{\tau}(n):=\min \left\{ i\geq 0:\tilde{S}_{i}=\tilde{L}_{n}\right\}
\end{equation*}%
and write%
\begin{equation*}
d_{n}=\sum_{k=0}^{n}\mathbf{E}\left[ \prod_{i=1}^{n}G_{i}\left(
F_{i,0}(0)\right) ;\tilde{\tau}(n)=k\right] .
\end{equation*}%
Recalling point 1) of Hypothesis A2 we conclude that, for any $i\geq 1$
\begin{equation*}
F_{i,0}^{\sigma }(0)=F_{i,i-1}^{\sigma }(F_{i-1,0}(0))\geq F_{i,i-1}^{\sigma
}(0)\geq \kappa ^{\sigma }.
\end{equation*}%
This estimate, point 2) of Hypothesis A2 and Lemma \ref{L_represent} imply
\begin{eqnarray*}
&&\mathbf{E}\left[ \prod_{i=1}^{n}G_{i}\left( F_{i,0}(0)\right) ;\tilde{\tau}%
(n)=k\right] \leq \mathbf{E}\left[ \prod_{i=1}^{n}G_{i}\left(
F_{i,0}^{\sigma }(0)\right) ;\tilde{\tau}(n)=k\right] \\
&&\quad \leq \mathbf{E}\left[ \left( \prod_{i=1}^{n}F_{i,0}^{\sigma
}(0)\right) ^{\gamma };\tilde{\tau}(n)=k\right] =\mathbf{E}\left[ \frac{1}{%
\left( \tilde{B}_{n}\right) ^{\sigma \gamma }};\tilde{\tau}(n)=k\right] .
\end{eqnarray*}%
Further,%
\begin{equation*}
\mathbf{E}\left[ \frac{1}{\left( \tilde{B}_{n}\right) ^{\sigma \gamma }};%
\tilde{\tau}(n)=k\right] \leq \mathbf{E}\left[ e^{\sigma \gamma \tilde{S}%
_{k}};\tilde{\tau}(n)=k\right] =\mathbf{E}\left[ e^{\sigma \gamma \tilde{S}%
_{k}};\tilde{\tau}(k)=k\right] \mathbf{P}\left( \tilde{L}_{n-k}\geq 0\right)
.
\end{equation*}%
Using Lemma \ref{Lappbasic} with $u(x)=e^{-\sigma \gamma x}$ we conclude
that, for any $\varepsilon >0$ there exists $m=m\left( \varepsilon \right) $
such that
\begin{eqnarray}
&&\sum_{k=m}^{n}\mathbf{E}\left[ \prod_{i=1}^{n}G_{i}\left(
F_{i,0}(0)\right) ;\tilde{\tau}(n)=k\right]  \notag \\
&&\quad \leq \sum_{k=m}^{n}\mathbf{E}\left[ e^{\sigma \gamma \tilde{S}_{k}};%
\tilde{\tau}(k)=k\right] \mathbf{P}\left( \tilde{L}_{n-k}\geq 0\right) \leq
\varepsilon \mathbf{P}\left( \tilde{L}_{n}\geq 0\right) .
\label{Estim_epsilon}
\end{eqnarray}

We now consider fixed $k\leq m$ and write
\begin{eqnarray*}
&&\mathbf{E}\left[ \prod_{i=1}^{n}G_{i}\left( F_{i,0}(0)\right) ;\tilde{\tau}%
(n)=k\right] \\
&&\qquad =\mathbf{E}\left[ \prod_{i=1}^{k}G_{i}\left( F_{i,0}(0)\right)
\prod_{j=k+1}^{n}G_{j}\left( F_{j,k}(F_{k,0}(0))\right) ;\tilde{\tau}(n)=k%
\right] \\
&&\qquad =\mathbf{E}\left[ \prod_{i=1}^{k}G_{i}\left( F_{i,0}(0)\right)
\Theta \left( n-k;F_{k,0}(0)\right) ;\tilde{\tau}(k)=k\right] ,
\end{eqnarray*}%
where%
\begin{equation*}
\Theta \left( n;s\right) :=\mathbf{E}\left[ \prod_{j=1}^{n}G_{j}\left(
F_{j,0}(s)\right) ;\tilde{L}_{n}\geq 0\right] .
\end{equation*}

Using the arguments applied to establish Lemma 2.7 in \cite{4h}, one may
check that, under the conditions of Theorem \ref{T_infvar0}
\begin{eqnarray*}
&&\sum_{j=1}^{\infty }\left( 1-G_{j}\left( F_{j,0}(s)\right) \right) \leq
\sum_{j=1}^{\infty }G_{j}^{\prime }(1)\left( 1-F_{j,0}(s)\right) \\
&&\quad\leq\sum_{j=1}^{\infty }G_{j}^{\prime }(1)\left( 1-F_{j,0}(0)\right)
\leq \sum_{j=1}^{\infty }G_{j}^{\prime }(1)e^{-\tilde{S}_{j}}<\infty \quad
\mathbf{\tilde{P}}^{+}-a.s.
\end{eqnarray*}

Hence it follows that,
\begin{equation*}
\xi _{n}(s):=\prod_{j=1}^{n}G_{j}\left( F_{j,0}(s)\right) \rightarrow \xi
_{\infty }(s):=\prod_{j=1}^{\infty }G_{j}\left( F_{j,0}(s)\right) >0
\end{equation*}%
$\mathbf{\tilde{P}}^{+}-$a.s. Since $\xi _{n}(s)\rightarrow \xi _{\infty }(s)
$ $\mathbf{\tilde{P}}^{+}-$a.s. as $n\rightarrow \infty $, it follows from
Lemma \ref{Lappbasic1} that, for each $s\in \lbrack 0,1)$
\begin{equation*}
\Theta \left( n;s\right) \sim \mathbf{\tilde{E}}^{+}\left[ \xi _{\infty }(s)%
\right] \mathbf{P}\left( \tilde{L}_{n}\geq 0\right) ,\,n\rightarrow \infty .
\end{equation*}

Applying the dominated convergence theorem gives on account of (\ref{Ladder2}%
) and properties of slowly varying functions%
\begin{eqnarray}
&&\lim_{n\rightarrow \infty }\mathbf{E}\left[ \prod_{i=1}^{k}G_{i}\left(
F_{i,0}(0)\right) \frac{\Theta \left( n-k;F_{k,0}(0)\right) }{\mathbf{P}%
\left( \tilde{L}_{n}\geq 0\right) };\tilde{\tau}(k)=k\right]  \notag \\
&&\quad =\mathbf{E}\left[ \prod_{i=1}^{k}G_{i}\left( F_{i,0}(0)\right)
\mathbf{\tilde{E}}^{+}\left[ \prod_{j=0}^{\infty }\hat{G}_{j}\left( \hat{F}%
_{j,0}(F_{k,0}(0))\right) \right] ;\tilde{\tau}(k)=k\right] ,
\label{Exact_limit}
\end{eqnarray}%
where $\hat{G}_{j},\hat{F}_{j,0}$ are independent copies of $G_{j},F_{j,0}.$

Combining (\ref{Exact_limit}) with (\ref{Estim_epsilon}) we get%
\begin{equation*}
\lim_{n\rightarrow \infty }\frac{1}{\mathbf{P}\left( \tilde{L}_{n}\geq
0\right) }\mathbf{E}\left[ \prod_{i=0}^{n-1}G_{i}\left( F_{i,n}(0)\right) %
\right] =\theta ,
\end{equation*}%
where%
\begin{equation*}
\theta :=\sum_{k=0}^{\infty }\mathbf{E}\left[ \prod_{i=1}^{k}G_{i}\left(
F_{i,0}(0)\right) \mathbf{\tilde{E}}^{+}\left[ \prod_{j=0}^{\infty }\hat{G}%
_{j}\left( \hat{F}_{j,0}(F_{k,0}(0))\right) \right] ;\tilde{\tau}(k)=k\right]
.
\end{equation*}

This proves Lemma \ref{L_infvarFrac}.

\textbf{Proof of Theorem \ref{T_infvar0}. }We know that
\begin{equation*}
d_{n}\sim \theta \frac{l_{3}(n)}{n^{\rho }}
\end{equation*}%
as $n\rightarrow \infty .$ This and a Tauberian theorem (see \cite{FE},
Chapter XIII.5, Theorem 5) imply
\begin{equation*}
D(s)=\sum_{n=1}^{\infty }d_{n}s^{n}\sim \theta \Gamma \left( 1-\rho \right)
\frac{l_{3}\left( 1/(1-s)\right) }{\left( 1-s\right) ^{1-\rho }}.
\end{equation*}%
Thus,%
\begin{equation*}
\mathcal{R}(s)=\frac{s\left( \mathcal{H}^{\ast }(s)+R_{1}\right) }{\left(
1-s\right) D\left( s\right) }\sim \frac{\mathcal{H}^{\ast }(1)+R_{1}}{\theta
\Gamma \left( 1-\rho \right) l_{3}\left( 1/(1-s)\right) \left( 1-s\right)
^{\rho }}
\end{equation*}%
as $s\uparrow 1.$ Since the sequence $\left\{ R_{n},n\geq 1\right\} $ is
monotone decreasing, it follows that (see \cite{FE}, Chapter XIII.5, Theorem
5)
\begin{equation*}
R_{n}\sim \frac{\mathcal{H}^{\ast }(1)+R_{1}}{\theta \Gamma \left( \rho
\right) \Gamma \left( 1-\rho \right) }\frac{n^{\rho -1}}{l_{3}\left(
n\right) }\quad \text{ as }n\rightarrow \infty .
\end{equation*}%
Theorem \ref{T_infvar0} is proved.

\end{document}